\documentclass[11pt, british]{amsart}

\usepackage{babel,eucal,url,amssymb,
enumerate,amscd,
}

\newtheorem{thm}{Theorem}[section]

\newtheorem{pro}[thm]{Proposition}
\newtheorem{co}[thm]{Corollary}

\newtheorem{defn}[thm]{Definition}

\newcommand{\Gtwo}{\ifmmode{{\rm G}_2}\else{${\rm G}_2$}\fi}

\newcommand{\g}{\mathfrak{g}}

\newcommand{\F}{\mathcal{F}}
\newcommand{\R}{\mathbb{R}}

\newcommand{\ee}{\end{equation}}
\newcommand{\be}[1]{\begin{equation}\label{#1}}

\begin{document}

\title[]
 {Totally umbilical  radical screen transversal half lightlike submanifolds of almost contact B-metric manifolds}


\author[G. Nakova]{Galia Nakova}

\address{
St. Cyril and St. Methodius University of Veliko Tarnovo, Bulgaria}
\email{gnakova@gmail.com}




\begin{abstract}
The present paper is a continuation of our previous work, where a class of half lightlike submanifolds of almost contact B-metric manifolds was introduced. We study curvature properties of totally and screen totally umbilical such submanifolds as well as of the corresponding semi-Riemannian submanifolds with respect to the associated B-metric.
\end{abstract}

\maketitle

\section{Introduction}\label{sec-1}
Many authors have studied half lightlike submanifolds of semi-Riemannian manifolds \cite{D-B,D-S,J}. In \cite{GN} we defined radical screen transversal half lightlike (RSTHL) submanifolds of almost contact B-metric manifolds and have examined the class of ascreen such submanifolds. This paper deals with curvature properties of totally umbilical and screen totally umbilical of the considered submanifolds. 
The main result concerning an ascreen RSTHL submanifold $(M,g)$  and its corresponding semi-Riemannian submanifold $(M,\widetilde g)$ is presented in the following theorem:
\begin{thm}\label{Theorem 4.6}
Let $(M,g)$ be a screen totally umbilical ascreen RSTHL submanifold of $\overline M(\overline \nu (p),\overline {\widetilde \nu }(p))$ and $\overline \nu \neq 0$. Then the following assertions are equivalent:\\
(i) $(M,g)$ is Ricci semi-symmetric; \quad (ii) $(M,\widetilde g)$ is Ricci semi-symmetric;\\
(iii) $(M,g)$ is $\overline \eta $-Einstein;\quad (iv) $(M,\widetilde g)$ is Einstein;\quad (v) $\overline \nu =4\mu ^2\gamma ^2$ on $M$.
\end{thm}
An example that confirms Theorem \ref{Theorem 4.6} is constructed.
\section{Preliminaries}\label{sec-2}
A $(2n+1)$-dimensional smooth manifold  $\overline M$ has an almost contact structure $(\overline \varphi,\overline \xi,\overline \eta)$ if it admits a $(1,1)$ tensor field $\overline \varphi $, a vector field 
$\overline \xi $ and a 1-form $\overline \eta $, satisfying the following conditions:
\begin{equation*}
\overline \varphi^2X=-X+\overline \eta(X)\overline \xi, \qquad \quad \overline \eta(\overline \xi)=1, \quad X\in T\overline M. 
\end{equation*}
If $(\overline M,\overline \varphi,\overline \xi,\overline \eta )$ is equipped with a semi-Riemannian metric $\overline g$, called {\it a B-metric} \cite{GaMGri}, such that
$\overline g(\overline \varphi X,\overline \varphi Y)=-\overline g(X,Y)+\overline \eta(X)\overline \eta(Y)$, then $(\overline M,\overline \varphi,\overline \xi,\overline \eta,\overline g)$ is called {\it an almost contact B-metric manifold}  \cite{GaMGri}.
Immediate consequences of the above conditions are:
\begin{equation*}
\overline \eta \circ \overline \varphi =0, \quad \overline \varphi \overline \xi =0, \quad {\rm rank}(\overline \varphi)=2n, \quad \overline \eta (X)=
\overline g(X,\overline \xi ), \quad \overline g(\overline \xi,\overline \xi )=1.
\end{equation*}
The tensor field $\overline {\widetilde g}$ of type $(0,2)$ given by 
$\overline {\widetilde g}(X,Y)=\overline g(X,\overline \varphi Y)+\overline \eta (X)\overline \eta (Y)$
is a B-metric, called {\it an associated metric} to $\overline g$. Both metrics $\overline g$ and 
$\overline {\widetilde g}$ are necessarily of signature $(n+1,n)$. 
Throughout this paper, for the orthonormal basis the signature of the metric $g$ will be of the form $(+\ldots + -\ldots -)$.\\
Let $\overline \nabla$ be the Levi-Civita connection of $\overline g$. A classification of the almost contact B-metric manifolds with respect to the tensor $F(X,Y,Z)=\overline g((\overline \nabla_X\overline \varphi)Y,Z)$ is given in \cite{GaMGri} and eleven basic classes  $\F_i$ $(i=1,2,\dots,11)$ are obtained. If $(\overline M,\overline \varphi,\overline \xi,\overline \eta,\overline g)$ belongs to $\F_i$ then it is called an {$\F_i$-{\it manifold}.
The special class $\F_0$ is determined by the condition $F(X,Y,Z)=0$ and in this class we have
$\overline \nabla \overline \varphi =\overline \nabla \overline \xi =\overline \nabla \overline \eta =\overline \nabla \overline g=\overline \nabla \overline {\widetilde g}=0$.
Let $\overline {\widetilde \nabla }$ be the Levi-Civita connection of $\overline {\widetilde g}$. 
In \cite{GaMGri} it is shown that 
the Levi-Civita connections $\overline \nabla $ and $\overline {\widetilde \nabla }$ of an $\F_0$-manifold coincide.
\par
Further, we briefly recall the main notions about half lightlike submanifolds of semi-Riemannian manifolds for which we refer to \cite{D-B, D-S}.
\par
A lightlike submanifold $(M,g)$ of codimension 2 of $(\overline M,\overline g)$ is called
{\it a half lightlike submanifold} if {\it the radical distribution}
${\rm Rad} (TM): p\in M \longrightarrow {\rm Rad} (T_pM)$ of $M$ has rank 1, where the {\it radical subspace} ${\rm Rad} (T_pM)$ is defined by 
\[\rm{Rad (T_pM)}=\{\xi _p \in T_pM: g(\xi _p,X_p)=0, \, \, \forall X_p\in T_pM\} .\]
\noindent
Also, with a half lightlike submanifold $M$ the following  two non-degenerate distributions are related:
{\it a screen distribution} $S(TM)$ and a 1-dimensional {\it screen transversal bundle} $S(TM^\bot)$ such that $S(TM^\bot)={\rm span}\{L\}$, $\overline g(L,L)=\epsilon =\pm 1$. 
For the tangent bundle $TM$ and the normal bundle $TM^\bot $  the following decompositions are valid:
\begin{equation}\label{2.1}
TM={\rm Rad} (TM)\bot S(TM) , \quad TM^\bot ={\rm Rad} (TM)\bot S(TM^\bot ),
\end{equation}
where the symbol $\bot $ denotes the orthogonal direct sum. It is well known \cite{D-B, D-S} that 
for any $\xi \in \Gamma ({\rm Rad} (TM))$ there exists a unique locally defined vector field $N$ satisfying
$
\overline g(N,\xi )=1, \quad \overline g(N,N)=\overline g(N,L)=\overline g(N,X)=0,  \, \forall X\in \Gamma (S(TM))
$.
The 1-dimensional vector bundle ${\rm ltr}(TM)={\rm span}\{N\}$  is called {\it the lightlike transversal bundle} of $M$ with respect to $S(TM)$. The {\it transversal vector  bundle} ${\rm tr}(TM)$ of $M$ is the complementary (but never orthogonal) vector bundle to $TM$ in $T\overline M$ such that
${\rm tr}(TM)=S(TM^\bot)\bot {\rm ltr}(TM)$.
Thus, for $T\overline M$ we have
\begin{equation}\label{2.2}
T\overline M=TM\oplus {\rm tr}(TM)=\{{\rm Rad} (TM)\oplus {\rm ltr}(TM)\}\bot S(TM)\bot S(TM^\bot),
\end{equation}
where $\oplus $ denotes a non-orthogonal direct sum.
Denote by $P$ the projection of $TM$ on $S(TM)$, from the first decomposition in  \eqref{2.1} for any
$X \in \Gamma (TM)$ we obtain $X=PX+\eta (X)\xi $, where $\eta $ is a differential 1-form on $M$ given by $\eta (X)=\overline g(X,N)$.
\par
The local Gauss-Weingarten formulas of $(M,g)$ and $S(TM)$ are given by
\begin{equation*}\label{2.3}
\begin{array}{l}
\overline \nabla _XY=\nabla _XY+B(X,Y)N+D(X,Y)L ,
\end{array}
\end{equation*}
\begin{equation*}\label{2.4}
\begin{array}{l}
\overline \nabla _XN=-A_NX+\tau (X)N+\rho (X)L , \quad \overline \nabla _XL=-A_LX+\phi (X)N ;
\end{array}
\end{equation*}
\begin{equation*}\label{2.6}
\nabla _XPY=\nabla ^*_XPY+C(X,PY)\xi , \quad \nabla _X\xi =-A^*_\xi X-\tau (X)\xi , \, \,
X, Y \in \Gamma (TM).
\end{equation*}
The induced connections $\nabla $ and 
$\nabla ^*$ on $TM$ and $S(TM)$, respectively, are linear connections; $A_N$, $A_L$ and
$A^*_\xi $ are the shape operators on $TM$ and $S(TM)$, respectively, and $\tau $, $\rho $, 
$\phi $ are 1-forms on $TM$; $B$ and $D$ are the local second fundamental forms of 
$M$, $C$ is the local second fundamental form of $S(TM)$ and  $B$, $C$, $D$ are related to their shape operators as follows:
\begin{equation*}
B(X,Y)=g(A^*_\xi X,Y), \, \overline g(A^*_\xi X,N)=0; \, \, C(X,PY)=g(A_N X,PY), \, \overline g(A_N X,N)=0;
\end{equation*}
\begin{equation*}
\epsilon D(X,PY)=g(A_L X,PY), \quad  \overline g(A_L X,N)=\epsilon \rho (X),  
\end{equation*}
\begin{equation*}
\epsilon D(X,Y)=g(A_L X,PY)-\phi (X)\eta (Y), \qquad \forall X, Y\in \Gamma (TM).
\end{equation*}
Since $\overline \nabla $ is torsion-free, $\nabla $ is also torsion-free. Therefore $B$ and $D$ are symmetric $F(M)$-bilinear forms on $\Gamma (TM)$ such that
$B(X,\xi )=0, \quad D(X,\xi )=-\phi (X), \quad \forall X\in \Gamma (TM)$.
In general, the induced connection $\nabla $ is not metric and satisfies
\begin{equation*}
(\nabla _Xg)(Y,Z)=B(X,Y)\eta (Z)+B(X,Z)\eta (Y).
\end{equation*}
The linear connection $\nabla ^*$ is not torsion-free but it is a metric connection on $S(TM )$,
$A^*_\xi $ and $A_N$ are $\Gamma (S(TM))$-valued, $A^*_\xi $ is self-adjoint with respect to $g$ and $A^*_\xi \xi =0$.
\par
In \cite{GN} we introduced {\it a Radical Screen Transversal  Half Lightlike (RSTHL) submanifold} of an almost contact B-metric manifold. In the remainder of this section we provide main notions and formulas for {\it ascreen} RSTHL, obtained in \cite{GN}. We recall:\\
A half lightlike submanifold $M$ of an almost  contact B-metric manifold $(\overline M,\overline \varphi,\overline \xi,\overline \eta,\overline g)$ is said to be  RSTHL \cite{GN} if $\overline \varphi ({\rm Rad} (TM))=S(TM^\bot )$. \\
A half lightlike submanifold $M$ of an almost  contact B-metric manifold $(\overline M,\overline \varphi,\overline \xi,\overline \eta,\overline g)$ is said to be ascreen  if $\overline \xi $  belongs to ${\rm Rad} (TM)\oplus {\rm ltr} (TM)$ \cite{J}. 
According to the notations in \cite{GN}, for an ascreen RSTHL submanifold $(M,g)$ of an almost contact B-metric manifold $(\overline M,\overline \varphi,\overline \xi,\overline \eta,\overline g)$ we have:
\begin{itemize}
\item $\overline \varphi \xi =\mu L$, \, \, $\overline \eta (\xi )=\mu $, $\forall \xi \in \Gamma ({\rm Rad} (TM))$, $\mu \neq 0$ is a smooth function on $M$;
\item $\overline g(L,L)=1$, \, \, $\overline \eta (L)=0$;
\item $\overline \xi =(1/2\mu )\xi +\mu N$, \, \, $\overline \varphi (S(TM))=S(TM)$, \, \, $\overline \eta (X)=\mu \eta (X), \, X\in \Gamma (TM)$;
\item $\overline \eta (N)=1/2\mu $, \, \, $\overline \varphi N=-(1/2\mu )L$, \, \, $\overline \varphi L=-(1/2\mu )\xi +\mu N$.  
\end{itemize}
Since on an almost contact B-metric manifold $\overline M$ there exist two B-metrics $\overline g$ and 
$\overline {\widetilde g}$, in \cite[Theorem 1.1]{GN} we consider two induced metrics $g$ and 
$\widetilde g$ on a submanifold $M$ of $\overline M$ by $\overline g$ and $\overline {\widetilde g}$, respectively, and we proved the following result: If $\overline M$ is a $(2n+1)$-dimensional almost contact B-metric manifold and $(M,g,S(TM),{\rm Rad} (TM))$ is an ascreen RSTHL submanifold of $(\overline M,\overline g)$, then 
$(M ,\widetilde g)$ is a semi-Riemannian submanifold of $(\overline M,\overline {\widetilde g})$ of codimension two. Moreover, the tangent bundle $TM$ of $(M,\widetilde g)$ is an orthogonal direct sum
with respect to $\widetilde g$ of the non-degenerate with respect to $\widetilde g$
distributions $S(TM)$ and ${\rm Rad} (TM)$, ${\rm Rad} (TM)$ is spacelike and the signature of 
$\widetilde g$ on $S(TM)$ is $(n-1,n-1)$.
\par
From now on, $(M,g)$ and $(M,\widetilde g)$ are the submanifolds of $\overline M$ from \cite[Theorem 1.1]{GN} and $\overline M$ is an $\F_0$-manifold.  For the induced geometric objects on 
$(M,g)$ we have \cite{GN}:
\begin{equation}\label{2.7}
A_NX=-(1/2\mu ^2)A^*_\xi X, \quad A_LX=(1/\mu )\overline \varphi (A^*_\xi X); 
\end{equation}
\begin{equation}\label{2.8}
D(X,Y)=(1/\mu )B(X,\overline \varphi (PY)), \quad C(X,PY)=-(1/2\mu ^2)B(X,Y); 
\end{equation}
\begin{equation}\label{2.9}
\tau (X)=-X(\mu )/\mu , \quad \phi (X)=\rho (X)=0, \,\, \forall X, Y \in \Gamma (TM);
\end{equation}
\begin{equation}\label{2.10}
(\nabla ^*_X\overline \varphi )(PY)=0 , \, \, \forall X,Y\in \Gamma (TM);
\end{equation}
\begin{equation*}
\begin{array}{ll}
\text{The shape operators}\, A^*_\xi , A_N \, \text{and}\, A_L\, \text{commute with}\, \overline \varphi
 \, \text{on} \, S(TM)\,  \text{and} \\
 B(\overline \varphi X,\overline \varphi Y)=-B(X,Y), \, \forall X,Y\in S(TM). 
\end{array}
\end{equation*}
The Gauss-Weingarten formulas of the submanifold $(M,\widetilde g)$ of $(\overline M,\overline {\widetilde g})$ are \cite{GN}: 
\begin{equation*}
\begin{array}{lll}
\overline {\widetilde \nabla }_XY=\widetilde \nabla _XY+h_1(X,Y)N_1+h_2(X,Y)N_2 , \\
\overline {\widetilde \nabla }_XN_1=-\widetilde A_{N_1}X , \quad \overline {\widetilde \nabla }_XN_2=-\widetilde A_{N_2}X , \, \, \forall X,Y\in \Gamma (TM),
\end{array}
\end{equation*}
where: $\widetilde \nabla $ is the Levi-Civita connection of $\widetilde g$; $N_1=\overline \xi -L$, $N_2=2\overline \xi-2\mu N-L$ are normal vector fields to $(M,\widetilde g)$ satisfying $\overline {\widetilde g}(N_1,N_1)=-
\overline {\widetilde g}(N_2,N_2)=1$, $\overline {\widetilde g}(N_1,N_2)=0$;
$\widetilde h(X,Y)=h_1(X,Y)N_1+h_2(X,Y)N_2$ is the second fundamental form of $(M,\widetilde g)$; 
$\widetilde A_{N_1}$ and $\widetilde A_{N_2}$ are  the shape operators with respect to $N_1$ and 
$N_2$, respectively; $h_1(X,Y)=\widetilde g(\widetilde A_{N_1}X,Y)$, \, $h_2(X,Y)=-\widetilde g(\widetilde A_{N_2}X,Y)$, are bilinear symmetric forms. 
The induced geometric objects on $(M,g)$ and $(M,\widetilde g)$ are related as follows \cite{GN}:
\begin{equation}\label{2.11}
\widetilde \nabla _XY=\nabla _XY+(1/\mu ^2)\left((1/2)B(X,Y)+B(X,\overline \varphi (PY))\right)\xi ;
\end{equation}
\begin{equation}\label{2.12}
h_1(X,Y)=(1/\mu )B(X,Y), \, \, h_2(X,Y)=-(1/\mu )(B(X,Y)+B(X,\overline \varphi (PY)));
\end{equation}
\begin{equation*}
\widetilde A_{N_1}X=-(1/\mu )\overline \varphi (A^*_\xi X), \quad
\widetilde A_{N_2}X=(1/\mu )(A^*_\xi X-\overline \varphi (A^*_\xi X)).
\end{equation*}
\section{Totally umbilical submanifolds $(M,g)$ and $(M,\widetilde g)$ of an $\F_0$-manifold}\label{sec-3}
\begin{defn}\label{Definition 3.1}\cite{D-B}
A half lightlike submanifold $(M,g)$ of $\overline M$ is called {\it totally umbilical} in $\overline M$ if there is a smooth vector field $H$ on ${\rm tr}(TM)$ on any coordinate neighborhood $U\subset \overline M$ such that for  the second fundamental form $h(X,Y)=B(X,Y)N+D(X,Y)L$ of $M$ the equality $h(X,Y)=Hg(X,Y)$ holds for any $X,U\in \Gamma (TM)$. In case $H=0$ ($H\neq 0$) on $U$, $M$ is {\it totally geodesic} ({\it proper totally umbilical}).
\end{defn}
It is clear that $M$ is totally umbilical if and only if $B(X,Y)=\beta g(X,Y)$ and $D(X,Y)=\delta g(X,Y)$ for any $X,Y\in \Gamma (TM)$, where $\beta $ and $\delta $ are smooth functions.
\begin{defn}\label{Definition 3.2}\cite{D-B}
A screen distribution $S(TM)$ of a half lightlike submanifold $(M,g)$ of $\overline M$ is called {\it totally umbilical} in $M$ if there is a smooth function $\gamma $on any coordinate neighborhood $U\subset M$ such that $C(X,PY)=\gamma g(X,Y)$ for any $X,Y\in \Gamma (TM)$. In case $\gamma =0$ ($\gamma \neq 0$) on $U$, $S(TM)$ is {\it totally geodesic} ({\it proper totally umbilical}) in $M$.
\end{defn}
A half lightlike submanifold $M$ with a totally umbilical $S(TM)$ we will call briefly  {\it screen totally umbilical}. By virtue of \eqref{2.8}, \eqref{2.12} we establish the truthfulness of
\begin{pro}\label{Proposition 3.3}
For the submanifolds $(M,g)$ and $(M,\widetilde g)$ of $\overline M\in \F_0$ we have:
\begin{enumerate}
\item[(i)] $(M,g)$ is totally geodesic if and only if $(M,\widetilde g)$ is totally geodesic.
\item[(ii)] If $(M,g)$ is totally umbilical, then  $(M,g)$, $S(TM)$ and $(M,\widetilde g)$ are totally geodesic.
\item[(iii)] If $(M,\widetilde g)$ is totally umbilical, then  $(M,\widetilde g)$, $S(TM)$ and  $(M,g)$ are totally geodesic.
\end{enumerate}
\end{pro}
The Ricci tensor of a lightlike submanifold (defined by $Ric(Y,Z)={\rm trace}\{X\longrightarrow R(X,Y,Z)\}$, where $R$ is the curvature tensor of $\nabla $)  is not symmetric in general because the induced connection $\nabla $ is not metric. It is known \cite{D-B, D-S} that a necessary and sufficient condition the Ricci tensor on a half lightlike submanifold $M$ to be symmetric is the 1-form $\tau $ to be closed, i.e. 
$d\tau =0$. Note that in our case $Ric$ of $(M,g)$ is symmetric (see the first equality in \eqref{2.9}).
\begin{pro}\label{Proposition 3.4}
The curvature tensors $R$, $\widetilde R$ and the Ricci tensors $Ric$, $\widetilde {Ric}$ \\of the submanifolds $(M,g)$, $(M,\widetilde g)$, respectively, of $\overline M\in \F_0$ are related as follows:
\begin{align}\label{13}
\begin{aligned}
\widetilde R(X,Y,Z)=R(X,Y,Z)+[B(Y,Z)+2B(Y,\overline \varphi (PZ))]A_NX\\
-[B(X,Z)+2B(X,\overline \varphi (PZ))]A_NY\\
+\left\{(1/2\mu ^2)[(\nabla _XB)(Y,Z)-(\nabla _YB)(X,Z)+\tau (X)B(Y,Z)-\tau (Y)B(X,Z)]\right. \\
+(1/\mu ^2)\left[\tau (X)B(Y,\overline \varphi (PZ))-\tau (Y)B(X,\overline \varphi (PZ))\right. \\
\left.\left. +(\nabla _XB)(Y,\overline \varphi (PZ))-(\nabla _YB)(X,\overline \varphi (PZ))\right]\right\}\xi ,
\end{aligned}
\end{align}
\begin{align}\label{14}
\begin{aligned}
\widetilde Ric(Y,Z)=Ric(Y,Z)+[B(Y,Z)+2B(Y,\overline \varphi (PZ))]{\rm trace} A_N-B(A_NY,Z)\\
-2B(A_NY,\overline \varphi (PZ))+(1/2\mu ^2)[(\nabla _\xi B)(Y,Z)-(\nabla _YB)(\xi ,Z)+\tau (\xi )B(Y,Z)]\\
+(1/\mu ^2)[(\nabla _\xi B)(Y,\overline \varphi (PZ))-(\nabla _YB)(\xi ,\overline \varphi (PZ))+
\tau (\xi )B(Y,\overline \varphi (PZ))].
\end{aligned}
\end{align}
\end{pro}
As an immediate consequence of Proposition~\ref{Proposition 3.3} and Proposition~\ref{Proposition 3.4} we state
\begin{co}\label{Corollary 3.5}
If $(M,g)$ or $(M,\widetilde g)$ is totally umbilical,, then $R=\widetilde R$, $Ric=\widetilde Ric$.
\end{co}
\section{Screen totally umbilical ascreen RSTHL submanifolds of $\F_0$-manifolds of constant totally real sectional curvatures}\label{sec-4}
The condition $\overline \nabla =\overline {\widetilde \nabla }$ on an $\F_0$-manifold $\overline M$ implies that the curvature tensors $\overline R$ and $\overline {\widetilde R}$ of type $(1,3)$ of 
$\overline \nabla$ and $\overline {\widetilde \nabla }$, respectively, coincide. The corresponding curvature tensors $\overline R$, $\overline {\widetilde R}$ of type $(0,4)$ are given by 
$\overline R(X,Y,Z,W)=\overline g(\overline R(X,Y,Z),W)$, $\overline {\widetilde R}(X,Y,Z,W)=\overline {\widetilde g}(\overline {\widetilde R}(X,Y,Z),W)$ and in \cite{N-G} it was shown that $\overline {\widetilde R}(X,Y,Z,W)=\overline R(X,Y,Z,\overline \varphi W)=\overline R(X,Y,\overline \varphi Z,W)$. For every non-degenerate with respect to $\overline g$ section $\alpha ={\rm span}\{x,y\}$  in $T_p\overline M$ the following two sectional curvatures are defined: 
\[
\overline \nu (\alpha ;p)=\frac{\overline R(x,y,y,x)}{\overline g(x,x)\overline g(y,y)-(\overline g(x,y))^2} ,\quad
\overline {\widetilde \nu }(\alpha ;p)=\frac{\overline {\widetilde R}(x,y,y,x)}{\overline g(x,x)\overline g(y,y)-(\overline g(x,y))^2} .
\]
A section $\alpha $ is said to be {\it totally real} with respect to $\overline g$ if $\overline \varphi \alpha $ 
is orthogonal to $\alpha $ with respect to $\overline g$ and the curvatures $\overline \nu (\alpha )$ and
$\overline {\widetilde \nu }(\alpha )$ are called {\it totally real sectional curvatures}.
\begin{thm}\label{Theorem 4.1}\cite{N-G}
Let $(\overline M,\overline\varphi ,\overline \xi ,\overline \eta ,\overline g)$ ($\rm {dim} \overline M\geq 5$) be an $\F_0$-manifold and $\alpha $ be a non-degenerate totally real orthogonal to $\overline \xi $ section in $\overline M$ (with respect to $\overline g$). $\overline M$ is of pointwise constant  totally real sectional curvatures $\overline \nu (p)$ and $\overline {\widetilde \nu }(p)$ of $\alpha $ if and only if
\begin{equation*}
\overline R(X,Y,Z,W)=\overline \nu [\pi _1(\overline \varphi X,\overline \varphi Y,\overline \varphi Z,\overline \varphi W)-\pi _2(X,Y,Z,W)]+\overline {\widetilde \nu }\pi _3(\overline \varphi X,\overline \varphi Y,\overline \varphi Z,\overline \varphi W). 
\end{equation*}
Both functions $\overline \nu (p)$ and $\overline {\widetilde \nu }(p)$ are constants if $\overline M$ is
connected and $\rm{dim} \overline M\geq 7$.
\end{thm}
\noindent
The tensors $\pi _1$, $\pi _2$ and $\pi _3$ are given by: $\pi _1(X,Y,Z,W)=\overline g(Y,Z)\overline g(X,W)-$\\
$\overline g(X,Z)\overline g(Y,W) , \,\ 
\pi _2(X,Y,Z,W)=\pi _1(X,Y,\overline \varphi Z,\overline \varphi W)$,\, $\pi _3(X,Y,Z,W)=$\\
$-\overline g(Y,Z)\overline g(X,\overline \varphi W)+\overline g(X,Z)\overline g(Y,\overline \varphi W)-\overline g(Y,\overline \varphi Z)\overline g(X,W)+\overline g(X,\overline \varphi Z)\overline g(Y,W)$.\\
By $\overline M(\overline \nu (p),\overline {\widetilde \nu }(p))$  we denote an $\F_0$-manifold of pointwise constant totally real sectional curvatures of a non-degenerate totally real orthogonal to 
$\overline \xi $ section  with respect to $\overline g$. Everywhere in the following assertions 
$\rm {dim} \,\overline M(\overline \nu (p),\overline {\widetilde \nu }(p))=2n+1$ and $n\geq 2$.
\par
The curvature tensors $\overline R$ and $R$ of a semi-Riemannian manifold $\overline M$ and its half lightlike submanifold $M$ are related as follows \cite{D-S}:
\begin{equation*}
\begin{array}{lll}
\overline R(X,Y,Z)=R(X,Y,Z)+B(X,Z)A_NY-B(Y,Z)A_NX+D(X,Z)A_LY\\
-D(Y,Z)A_LX+\left\{(\nabla _XB)(Y,Z)-(\nabla _YB)(X,Z)+\tau (X)B(Y,Z)-\tau (Y)B(X,Z)\right.\\
\left.+\phi (X)D(Y,Z)-\phi (Y)D(X,Z)\right\}N\\
+\left\{(\nabla _XD)(Y,Z)-(\nabla _YD)(X,Z)+\rho (X)B(Y,Z)-\rho (Y)B(X,Z)\right\}L,  
\end{array}
\end{equation*}
for any $X,Y,Z\in \Gamma (TM)$. Now, by using the latter equality, \eqref{2.7}, \eqref{2.8}, \eqref{2.9} and Theorem \ref{Theorem 4.1}, we obtain the following proposition:
\begin{pro}\label{Proposition 4.2}
The curvature tensor $R$ of the submanifold $(M,g)$ of $\overline M(\overline \nu (p),\overline {\widetilde \nu }(p))$ is given by
\begin{align}\label{15}
\begin{aligned}
R(X,Y,Z)=-B(X,Z)A_NY+2B(X,\overline \varphi (PZ))\overline \varphi (A_NY)+B(Y,Z)A_NX\\
-2B(Y,\overline \varphi (PZ))\overline \varphi (A_NX)-[\overline \nu \overline g(\overline \varphi Y,\overline \varphi Z)+\overline {\widetilde \nu }g(Y,\overline \varphi Z)]PX\\
+[\overline \nu \overline g(\overline \varphi X,\overline \varphi Z)+\overline {\widetilde \nu }g(X,\overline \varphi Z)]PY-
[\overline \nu g(Y,\overline \varphi Z)-\overline {\widetilde \nu }\overline g(\overline \varphi Y,\overline \varphi Z)]\overline \varphi (PX)\\
+[\overline \nu g(X,\overline \varphi Z)-\overline {\widetilde \nu }\overline g(\overline \varphi X,\overline \varphi Z)]\overline \varphi (PY)+(1/2)\left\{\overline \nu \left[g(Y,Z)\eta (X)
\right.\right.\\
\left.\left.-g(X,Z)\eta (Y)\right]-\overline {\widetilde \nu }[g(Y,\overline \varphi Z)
\eta (X)-g(X,\overline \varphi Z)\eta (Y)]\right\}\xi
\end{aligned}
\end{align}
and the next equality is fulfilled
\begin{equation}\label{16}
(\nabla _XB)(Y,Z)-(\nabla _YB)(X,Z)+\tau (X)B(Y,Z)-\tau (Y)B(X,Z)=
\end{equation}
\begin{equation*}
\mu ^2\{\overline \nu \left[g(X,Z)\eta (Y)
-g(Y,Z)\eta (X)\right]-\overline {\widetilde \nu }[g(X,\overline \varphi Z)
\eta (Y)-g(Y,\overline \varphi Z)\eta (X)]\}.
\end{equation*}
\end{pro}
From Proposition \ref{Proposition 3.3}, \eqref{15} and \eqref{16} it follows the  truthfulness of the next
\begin{co}\label{Corollary 4.3}
If the submanifold $(M,g)$ of $\overline M(\overline \nu (p),\overline {\widetilde \nu }(p))$ is totally umbilical, then $(M,g)$ and $\overline M$ are flat.
\end{co}
Let us assume that $(M,g)$ is a screen totally umbilical ascreen RSTHL submanifold of 
$\overline M\in \F_0$. Then according to Definition \ref{Definition 3.2}, \eqref{2.7} and \eqref{2.8} we have
\begin{equation}\label{17}
A_NX=\gamma PX , \quad B(X,Y)=-2\mu ^2\gamma g(X,Y), \quad \forall X,Y\in \Gamma (TM).
\end{equation}
Further we  deal with such submanifold of $\overline M(\overline \nu (p),\overline {\widetilde \nu }(p))$.   
\begin{thm}\label{Theorem 4.4}
Let $(M,g)$ be a screen totally umbilical ascreen RSTHL submanifold of $\overline M(\overline \nu (p),\overline {\widetilde \nu }(p))$. Then $\overline {\widetilde \nu }=0$ and the function $\gamma $ from Definition \ref{Definition 3.2} satisfies the partial differential equations
\begin{equation}\label{18}
\overline \nu +2\tau (\xi )\gamma -2\xi (\gamma )-4\mu ^2\gamma ^2=0, \quad PX(\mu \gamma )=0, \, \, \forall X\in \Gamma (TM).
\end{equation}
Moreover, the curvature tensor $R$ and the Ricci tensor $Ric$ of $(M,g)$ are given by
\begin{align}\label{19}
\begin{aligned}
R(X,Y,Z)=[(\overline \nu -2\mu ^2\gamma ^2)g(Y,Z)-\overline \nu \overline \eta (Y)\overline \eta (Z)]PX\\
-[(\overline \nu -2\mu ^2\gamma ^2)g(X,Z)-\overline \nu \overline \eta (X)\overline \eta (Z)]PY+
(4\mu ^2\gamma ^2-\overline \nu )g(Y,\overline \varphi Z)\overline \varphi (PX)\\
-(4\mu ^2\gamma ^2-\overline \nu )g(X,\overline \varphi Z)\overline \varphi (PY)+
(\overline \nu /2)[g(Y,Z)\eta (X)-g(X,Z)\eta (Y)]\xi ,
\end{aligned}
\end{align}
\begin{align}\label{20}
\begin{aligned}
Ric(Y,Z)=[((4n-7)/2)\overline \nu -2(2n-5)\mu ^2\gamma ^2]g(Y,Z)-2(n-1)\overline \nu \overline \eta (Y)\overline \eta (Z)
\end{aligned}
\end{align}
and the function $\mu \gamma $ is a constant on $M$ if and only if $\overline \nu =4\mu ^2\gamma ^2$ on $M$.
\end{thm}
By direct calculations, using \eqref{13}, \eqref{14}, \eqref{17}, \eqref{19} and \eqref{20} we get 
\begin{pro}\label{Proposition 4.5}
Let $(M,g)$ be a screen totally umbilical ascreen RSTHL submanifold of $\overline M(\overline \nu (p),\overline {\widetilde \nu }(p))$. Then the curvature tensor 
$\widetilde R$ and the Ricci tensor $\widetilde {Ric}$ of $(M,\widetilde g)$ are given by
\begin{align}\label{21}
\begin{aligned}
\widetilde {R}(X,Y,Z)=[(\overline \nu -4\mu ^2\gamma ^2)g(Y,Z)-4\mu ^2\gamma ^2g(Y,\overline \varphi Z)-\overline \nu \overline \eta (Y)\overline \eta (Z)]PX\\
-[(\overline \nu -4\mu ^2\gamma ^2)g(X,Z)-4\mu ^2\gamma ^2g(X,\overline \varphi Z)-\overline \nu \overline \eta (X)\overline \eta (Z)]PY\\
-(\overline \nu -4\mu ^2\gamma ^2)g(Y,\overline \varphi Z)\overline \varphi (PX)+
(\overline \nu -4\mu ^2\gamma ^2)g(X,\overline \varphi Z)\overline \varphi (PY)\\
+\overline \nu[g(X,\overline \varphi Z)\eta (Y)-g(Y,\overline \varphi Z)\eta (X)]\xi ,
\end{aligned}
\end{align}
\begin{align}\label{22}
\begin{aligned}
\widetilde {Ric}(Y,Z)=2(n-2)(\overline \nu -4\mu ^2\gamma ^2)g(Y,Z)\\
-[\overline \nu+4(2n-3)\mu ^2\gamma ^2]g(Y,\overline \varphi Z)-2(n-1)\overline \eta (Y)\overline \eta (Z), \, X,Y\in \Gamma (TM).
\end{aligned}
\end{align}
\end{pro}
A half lightlike submanifold $M$ is said to be {\it Ricci semi-symmetric} \cite{D-S} if 
\begin{equation*}
\begin{array}{ll}
({\mathcal R}(X,Y)\cdot Ric)(X_1,X_2)
=-Ric(R(X,Y,X_1),X_2)-Ric(X_1,R(X,Y,X_2))=0.  
\end{array}
\end{equation*}
We recall that $(M,g)$ is {\it $\overline \eta $-Einstein} if $Ric=kg+c\overline \eta \otimes \overline \eta $, where $k$ and $c$ are constants.
{\bf Proof of Theorem \ref{Theorem 4.6}.}
(i)$\Leftrightarrow $(v): By virtue of \eqref{19} and \eqref{20} we find
\begin{align}\label{23}
\begin{aligned}
({\mathcal R}(X,Y)\cdot Ric)(X_1,X_2)=(2n-5)\overline \nu (\overline \nu /2-2\mu ^2\gamma ^2)\left[g(X,X_2)\overline \eta (Y)\overline \eta (X_1)\right.\\
\left.-g(Y,X_2)\overline \eta (X)\overline \eta (X_1)+g(X,X_1)\overline \eta (Y)\overline \eta (X_2)-
g(Y,X_1)\overline \eta (X)\overline \eta (X_2)\right].
\end{aligned}
\end{align}
Assume that $(M,g)$ is Ricci semi-symmetric. Take $Y=PY$, $X_2=PX_2$, $X=X_1=\xi $  in \eqref{23}, we obtain
$-(2n-5)\overline \nu (\overline \nu /2-2\mu ^2\gamma ^2)\mu ^2g(PY,PX_2)=0$. Since $S(TM)$ is non-degenerate, we get $\overline \nu =4\mu ^2\gamma ^2$. The implication (v)$\Rightarrow $(i) follows from \eqref{23}.\\
(ii)$\Leftrightarrow $(v): With the help of \eqref{21} and \eqref{22} we compute
\begin{align}\label{24}
\begin{aligned}
(\widetilde {\mathcal R}(X,Y)\cdot \widetilde {Ric})(X_1,X_2)=(2n-3)\overline \nu (\overline \nu -4\mu ^2\gamma ^2)\left[g(X,\overline \varphi X_2)\overline \eta (Y)\overline \eta (X_1)\right.\\
\left.-g(Y,\overline \varphi X_2)\overline \eta (X)\overline \eta (X_1)+g(X,\overline \varphi X_1)\overline \eta (Y)\overline \eta (X_2)-g(Y,\overline \varphi X_1)\overline \eta (X)\overline \eta (X_2)\right]\\
-2(n-2)(\overline \nu -4\mu ^2\gamma ^2)\left\{4\mu ^2\gamma ^2\left[g(X,\overline \varphi X_1)
g(Y,X_2)-g(Y,\overline \varphi X_1)g(X,X_2)\right.\right.\\
\left.+g(X,\overline \varphi X_2)g(Y,X_1)-g(Y,\overline \varphi X_2)g(X,X_1)\right]+\overline \nu \left[
g(Y,X_1)\overline \eta (X)\overline \eta (X_2)\right.\\
\left.\left.-g(X,X_1)\overline \eta (Y)\overline \eta (X_2)+g(Y,X_2)\overline \eta (X)\overline \eta (X_1)-g(X,X_2)\overline \eta (Y)\overline \eta (X_1)\right]\right\}.
\end{aligned}
\end{align}
Let $(M,\widetilde g)$ be Ricci semi-symmetric. We put $Y=PY$ and $X=X_1=\xi $  in \eqref{24} and taking into account that  $S(TM)$ is non-degenerate, we have 
$-(2n-3)\overline \nu (\overline \nu -4\mu ^2\gamma ^2)\mu ^2\overline \varphi (PY)-2(n-2)\overline \nu (\overline \nu -4\mu ^2\gamma ^2)\mu ^2PY=0$. Since $PY$ and $\overline \varphi (PY)$ are linear independent, the latter equality implies (v). The implication (v)$\Rightarrow $(ii) is obvious.
\par
By using that $\mu \gamma $ is a constant  iff $\overline \nu =4\mu ^2\gamma ^2$ on $M$ (Theorem \ref{Theorem 4.4}), it is not hard to see that (iii) and (iv) are equivalent to (v). Thus, we complete the proof.
\par
Note that an Einstein semi-Riemannian manifold is always Ricci semi-symmetric, but the converse is not true. Also, an Einstein lightlike manifold is not Ricci semi-symmetric in general.\\
{\bf Example 4.7}
Consider the  Lie group $G^{\prime }$ and its real Lie algebra ${\mathfrak {g}}^{\prime }$ from
\cite[Example 5.2]{GN2}.
Let $\{X_1,X_2,X_3,X_4\}$ be a global basis of left invariant vector fields on $G^{\prime }$. In \cite{GN2} we defined an almost complex structure $J$ and a left invariant metric $g^{\prime }$ on 
$G^{\prime }$ by $JX_1=X_3, \quad JX_2=X_4, \quad JX_3=-X_1, \quad JX_4=-X_2;$
\begin{equation*}
g^\prime (X_i,X_i)=-g^\prime (X_{i+1},X_{i+1})=1, \, (i=1,2); \, 
g^\prime (X_i,X_j)=0; \, i\neq j; \, (i,j=1,2,3,4).
\end{equation*}
The non-zero commutators of the basis vector fields and the non-zero components of the Levi-Civita connection $\nabla ^{\prime }$ of $g^{\prime }$ are: 
\begin{equation*}
\left[X_1,X_2\right]=-[X_3,X_4]=-2X_4\, ;\qquad \left[X_1,X_4\right]=-[X_2,X_3]=2X_2 .
\end{equation*}
\begin{equation*}
\begin{array}{ll}
\nabla ^{\prime }_{X_2}X_1=-\nabla ^{\prime }_{X_4}X_3=2X_4 \, ;  \quad 
\nabla ^{\prime }_{X_2}X_2=-\nabla ^{\prime }_{X_4}X_4=-2X_3 \, ; 
\\
\nabla ^{\prime }_{X_2}X_3=\nabla ^{\prime }_{X_4}X_1=-2X_2 \, ;  \quad 
\nabla ^{\prime }_{X_2}X_4=\nabla ^{\prime }_{X_4}X_2=2X_1 .
\end{array}
\end{equation*}
In \cite{GN2} we showed that $F^{\prime }(X^{\prime },Y^{\prime },Z^{\prime })=g^{\prime }((\nabla ^{\prime }_{X^{\prime }}J)Y^{\prime },Z^{\prime })=0$, $\forall X^{\prime },Y^{\prime },Z^{\prime }\in {\mathfrak {g}}^{\prime }$, i.e. $(G^{\prime },J,g^{\prime })$ is a K\"ahler-Norden manifold.
\par
Now, we denote by ${\mathfrak {g}}^{''}=\rm {span}\{E\}$ the real Lie algebra of the 1-dimensional Lie group ${\R}$ and define a left invariant metric $g^{''}$ on ${\R}$ by 
$g^{''}(E,E)=1$. Then for the Levi-Civita connection $\nabla ^{''}$ of $g^{''}$ the equality
$\nabla ^{''}_EE=0$ holds.
\par
Let us consider the 5-dimensional Lie group $\overline G=G^\prime \times {\R}$ whose  real Lie algebra  is $\overline {\mathfrak {g}}={\mathfrak {g}}^\prime \oplus {\mathfrak {g}}^{''}$. We define an almost contact structure $(\overline \varphi,\overline \xi,\overline \eta)$ and a left invariant B-metric 
$\overline g$ on $\overline G$ in the following way:
\begin{equation*}
\begin{array}{ll}
\overline \varphi (X^\prime ,aE)=(JX^\prime ,0), \quad \overline \xi =(0,E), \quad
\overline \eta (X^\prime ,aE)=a, \quad X^\prime \in {\mathfrak {g}}^{\prime }, \quad a,b\in \R.\\
\overline g((X^\prime ,aE),(Y^\prime ,bE))=g^\prime (X^\prime ,Y^\prime )+g^{''}(aE,bE)
\end{array}
\end{equation*}
Thus, $(\overline G,\overline \varphi,\overline \xi,\overline \eta ,\overline g)$ is a 5-dimensional almost contact B-metric manifold. It is known \cite{GN3} that for the commutator of any two vector fields 
$(X^\prime ,aE)$, $(Y^\prime ,bE)$ on $\overline G$ we have 
$[(X^\prime ,aE),(Y^\prime ,bE)]=[X^\prime ,Y^\prime ]+[aE,bE]=[X^\prime ,Y^\prime ]$ and 
the Levi-Civita connection $\overline \nabla $ of $\overline g$, $\nabla ^\prime $ and $\nabla ^{''}$ are related as follows:
\begin{equation*}
\overline g(\overline \nabla _{(X^\prime ,aE)}(Y^\prime ,bE),(Z^\prime ,cE))=g^{\prime }(\nabla ^{\prime }_{X^\prime }Y^\prime ,Z^\prime )+g^{''}(\nabla ^{''}_{aE}bE,cE)=g^{\prime }(\nabla ^{\prime }_{X^\prime }Y^\prime ,Z^\prime ).
\end{equation*}
By using the above equality it is easy to see that $\overline \nabla _{(X^\prime ,aE)}(Y^\prime ,bE)=
\nabla ^{\prime }_{X^\prime }Y^\prime $. Then we obtain $\overline F((X^\prime ,aE),(Y^\prime ,bE),(Z^\prime ,cE))=\overline g((\overline \nabla _{(X^\prime ,aE)}\overline \varphi )(Y^\prime ,bE),(Z^\prime ,cE))=0$, i.e. $(\overline G,\overline \varphi,\overline \xi,\overline \eta ,\overline g)$ is an ${\F}_0$-manifold. For the curvature tensor of type $(0,4)$ $\overline R$ of $\overline G$ we find
$\overline R=4(\pi _1\circ \overline \varphi-\pi _2)$. From Theorem \ref{Theorem 4.1} it follows that 
$\overline G$ is of constant totally real sectional curvatures $\overline \nu =4$ and $\overline {\widetilde \nu }=0$.\\
The subspace ${\g}$ of $\overline {\g}$ spanned by $\{E_1=(X_2,0),E_2=(X_4,0),\xi =(-\mu X_3,\mu E)\}$, $\mu \in{\R}, \mu \neq 0$ is a Lie subalgebra of $\overline {\g}$. Hence, the corresponding  to ${\g}$ Lie subgroup $G$ of $\overline G$ is a 3-dimensional submanifold of $\overline G$. The induced metric $g$ on $G$ by $\overline g$ is degenerate and ${\rm Rad} ({\g})={\rm span} \{\xi \}$. Thus 
$(G,g)$ is a half lightlike submanifold of $(\overline G,\overline g)$. We take 
the screen distribution $S({\g})$, the lightlike transversal bundle ${\rm ltr} ({\g})$ and the screen transversal bundle $S({\g}^\bot )$ of $(G,g)$ as follows: $S({\g})={\rm span}\{E_1,E_2\}$, 
${\rm ltr} ({\g})={\rm span}\{N=(1/2\mu )(X_3,E)\}$, $S({\g}^\bot )={\rm span}\{L=(X_1,0)\}$. Since $\overline \varphi \xi =\mu L$ and $\overline \xi =(1/2\mu )\xi +\mu N$, it follows that $(G,g)$ is an ascreen RSTHL submanifold. We directly check that the induced metric $\widetilde g$ on $G$ by 
$\overline {\widetilde g}$ is non-degenerate. We find $\overline \nabla _{E_1}N=-(1/\mu )E_1$,
$\overline \nabla _{E_2}N=-(1/\mu )E_2$, $\overline \nabla _\xi N=0$. The last three equalities imply $A_NX=(1/\mu )X$ for any vector field $X\in {\g}$. Hence, $C(X,PY)=(1/\mu )g(X,Y)$, 
which means that $(G,g)$ is screen totally umbilical and $\gamma =1/\mu $ and $\overline \nu =4\mu ^2\gamma ^2$. 
\par
The non-zero components of the linear connection $\nabla $ of $g$ are: $\nabla _{E_1}\xi =2\mu E_1$,
$\nabla _{E_2}\xi =2\mu E_2$, $\nabla _{E_1}E_1=-\nabla _{E_2}E_2=(1/\mu )\xi $.
Let $Y,Z\in {\g}$ and $Y=Y^1E_1+Y^2E_2+Y^3\xi $, $Z=Z^1E_1+Z^2E_2+Z^3\xi $. Then we get
$R(E_1,Y,Z)=-2(Y^2Z^2+2\mu ^2Y^3Z^3)E_1-2Y^2Z^1E_2-2Y^3Z^1\xi $, $R(E_2,Y,Z)=2Y^1Z^2E_1+2(Y^1Z^1-2\mu ^2Y^3Z^3)E_2+2Y^3Z^2\xi $, $R(\xi ,Y,Z)=4\mu ^2Y^1Z^3E_1+4\mu ^2Y^2Z^3E_2+2(Y^1Z^1-Y^2Z^2)\xi $. Now, for the Ricci tensor of $(G,g)$ we find $Ric(Y,Z)=4g(Y,Z)-8\overline \eta (Y)\overline \eta (Z)$.
From the latter equality it follows that $(G,g)$ is $\overline \eta $-Einstein. Analogously we obtain that the Ricci tensor $\widetilde {Ric}$ of $(G,\widetilde g)$ is given by $\widetilde {Ric}=-8\widetilde g(Y,Z)$, i.e. $(G,\widetilde g)$ is Einstein. Directly we check that both submanifolds $(G,g)$ and $(G,\widetilde g)$ of $\overline G$ are Ricci  semi-symmetric.

\end{document}